# Influence of the Forward Difference Scheme for the Time Derivative on the Stability of Wave Equation Numerical Solution

Aslam Abdullah
Faculty of Mechanical and Manufacturing Engineering
Universiti Tun Hussein Onn Malaysia
Batu Pahat, Malaysia
aslam@uthm.edu.my

*Abstract*— **Research on numerical stability of difference equations has been quite intensive in the past century. The choice of difference schemes for the derivative terms in these equations contributes to a wide range of the stability analysis issues - one of which is how a chosen scheme may directly or indirectly contribute to such stability. In the present paper, how far the forward difference scheme for the time derivative in the wave equation influences the stability of the equation numerical solution, is particularly investigated. The stability analysis of the corresponding difference equation involving four schemes, namely Lax's, central, forward, and rearward differences, were carried out, and the resulting stability criteria were compared. The results indicate that the instability of the solution of wave equation is not always due to the forward difference scheme for the time derivative. Rather, it is shown in this paper that the stability criterion is still possible when the spatial derivative is represented by an appropriate difference scheme. This sheds light on the degree of applicability of a difference scheme for a hyperbolic equation.**

I. INTRODUCTION

Finite difference is one of the numerical methods used in the discretization of partial differential equations such as Poisson's equation, wave equation, and Benjamin–Bona–Mahony (BBM) equation [1]-[3]. Following the discretization process, the stability analysis is necessary, in particular when it involves an explicit difference equation. This is to ensure whether or not the equation leads to a stable solution. Works done by, for instance, [4] and [5] serve as excellent materials for the analysis of stability.

*A. Discretization of the Hyperbolic Equation*

Without losing the generality, we begin with the first-order partial differential wave equation; in Cartesian coordinates and tensor notation following the Einstein convention, the equation is

$$\partial_t u + a \partial_x u = 0, \qquad (1)$$

where $a$ is a constant. The corresponding solution domain is covered by a grid shown in Figure 1.

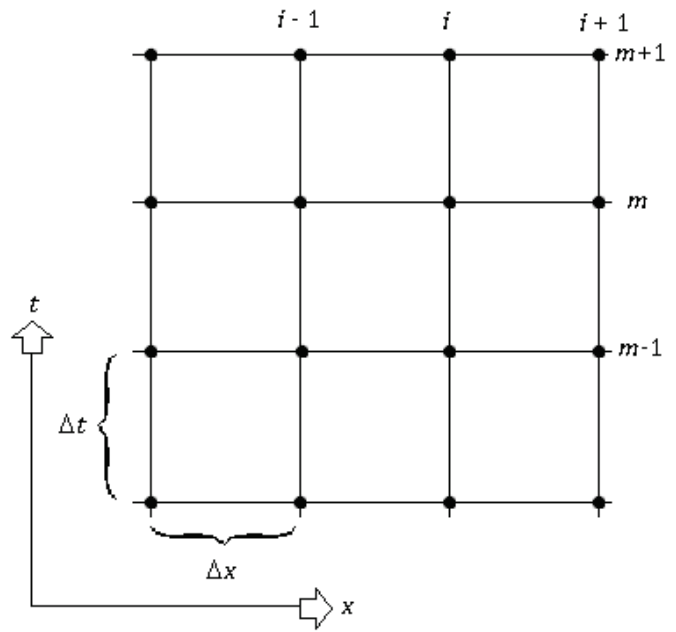

Fig. 1. Grid for the differencing of the first-order wave equation.

It is worth to note here that the grid expansion factor $r_e = 1$. The solution involving $r_e \neq 1$ has been considered in [6].

*B. Lax Method*

Replacing time and spatial derivatives with *Lax's difference* [7] and *central difference* schemes, respectively, we have, after some rearrangement

$$\frac{u_i^{m+1} - 1/2\left(u_{i+1}^m + u_{i-1}^m\right)}{\Delta t} = -a\frac{u_{i+1}^m - u_{i-1}^m}{2\Delta x} \qquad (2)$$

## II. ROUND-OFF ERROR

The round-off error is defined as

$$\varepsilon = N - E, \quad (3)$$

where $N$ and $E$ are finite accuracy numerical solution from a real computer and exact solution of difference equation, respectively. Note that the numerical solution $N$ satisfies the difference equation (2). Replacing $u$ in (2) with $N$ in (3), one has

$$\frac{E_i^{m+1} + \varepsilon_i^{m+1} - 1/2\left(E_{i+1}^m + \varepsilon_{i+1}^m + E_{i-1}^m + \varepsilon_{i-1}^m\right)}{\Delta t}$$
$$= -a\frac{E_{i+1}^m + \varepsilon_{i+1}^m - E_{i-1}^m - \varepsilon_{i-1}^m}{2\Delta x} \quad (4)$$

Similarly, replacing $u$ in (2) with $E$ gives

$$\frac{E_i^{m+1} - 1/2\left(E_{i+1}^m + E_{i-1}^m\right)}{\Delta t} = -a\frac{E_{i+1}^m - E_{i-1}^m}{2\Delta x}, \quad (5)$$

since $E$ clearly satisfies (2). It can be proven that $\varepsilon$ also satisfies the equation; subtracting (5) from (4), one has

$$\frac{\varepsilon_i^{m+1} - 1/2\left(\varepsilon_{i+1}^m + \varepsilon_{i-1}^m\right)}{\Delta t} = -a\frac{\varepsilon_{i+1}^m - \varepsilon_{i-1}^m}{2\Delta x} \quad (6)$$

## III. THE STABILITY

The solution of (2) is stable if and only if

$$\left|\frac{\varepsilon_i^{m+1}}{\varepsilon_i^m}\right| \leq 1 \quad (7)$$

We use a Fourier series to analytically represent the random variation of $\varepsilon$ with respect to space and time;

$$\varepsilon(x,t) = \sum_n e^{bt} e^{ik_n x}, \quad (8)$$

$$n = 1, 2, 3...$$

where $e^{bt}$ is the wave amplitude, $k_n$ is the wave number, and $b$ is a constant.

Note that the difference equation (2) is linear. Furthermore, it is satisfied by $\varepsilon$ as shown by (6). Thus, if the series (8) is substituted into (6), each term of the series and the series itself behave the same. This allows us to handle just one term of the series, namely

$$\varepsilon(x,t) = e^{bt} e^{ik_n x}, \quad (9)$$

in order to probe the stability of (2) without loss in generality. The stability of the difference equation depends on how $\varepsilon$ grows in time steps. Substituting (9) into (6) and (7) gives

$$\frac{e^{b(t+\Delta t)} e^{ik_n x} - 1/2\left(e^{bt} e^{ik_n(x+\Delta x)} + e^{bt} e^{ik_n(x-\Delta x)}\right)}{\Delta t}$$
$$= -a\frac{e^{bt} e^{ik_n(x+\Delta x)} - e^{bt} e^{ik_n(x-\Delta x)}}{2\Delta x} \quad (10)$$

and

$$\left|e^{b\Delta t}\right| \leq 1 \quad (11)$$

Dividing (10) by $e^{bt} e^{ik_n x}$;

$$\frac{e^{b\Delta t} - 1/2\left(e^{ik_n \Delta x} + e^{-ik_n \Delta x}\right)}{\Delta t} = -a\frac{e^{ik_n \Delta x} - e^{-ik_n \Delta x}}{2\Delta x} \quad (12)$$

Combining (11) and (12);

$$\left|\frac{1}{2}\left[-\frac{a\Delta t}{\Delta x}\left(e^{ik_n \Delta x} - e^{-ik_n \Delta x}\right) + \left(e^{ik_n \Delta x} + e^{-ik_n \Delta x}\right)\right]\right| \leq 1. \quad (13)$$

Let $e^{ik_n \Delta x} := \xi_1$, and $e^{-ik_n \Delta x} := \xi_2$. Rewriting (13), one has

$$\left|\frac{1}{2}\left[-\frac{a\Delta t}{\Delta x}(\xi_1 - \xi_2) + (\xi_1 + \xi_2)\right]\right| \leq 1 \quad (14)$$

Here we have two possible conditions which must hold simultaneously:

### A. Condition a

$$-\frac{1}{2}\left[\frac{a\Delta t}{\Delta x}(\xi_1 - \xi_2) - (\xi_1 + \xi_2)\right] \leq 1 \quad (15)$$

Setting $\xi_1 = -1$, and $\xi_2 = 1$, (15) becomes

$$\frac{a\Delta t}{\Delta x} \leq 1 \quad (16)$$

## B. Condition b

$$-\frac{1}{2}\left[\frac{a\Delta t}{\Delta x}(\xi_1 - \xi_2) - (\xi_1 + \xi_2)\right] \geq -1 \quad (17)$$

Setting $\xi_1 = 1$, and $\xi_2 = -1$, (17) becomes

$$\frac{a\Delta t}{\Delta x} \leq 1, \quad (18)$$

which is identical to (16) in *Condition a*.

The stability of the difference equation, (2), is therefore proven. Moreover, the stability requirement is given by (16) or (18) called the *Courant-Friedrichs-Lewy* (CFL) condition [8]-[10]. Any number given by the LHS term in (18) (i.e. $a\Delta t/\Delta x$) is called the *Courant number C* such that

$$C = \frac{a\Delta t}{\Delta x} \leq 1 \quad (19)$$

It is important to note that the Lax's difference scheme used to represent the time derivative in (1) is accurate to the first-order. The main reason a simpler first-order difference is not applicable is due to the instability of the resulting difference equations, an issue that will be discussed further in the next section. The more complex analysis of stability can be found in [11]-[15].

### IV. SIMPLE FIRST-ORDER DIFFERENCE FOR THE TIME DERIVATIVE

Replacing the time derivative with a simple first-order forward difference instead of *Lax's difference* schemes while keeping the *central difference* scheme to represent the spatial derivative in (1), we have, after some rearrangement

$$\frac{u_i^{m+1} - u_i^m}{\Delta t} = -a\frac{u_{i+1}^m - u_{i-1}^m}{2\Delta x} \quad (20)$$

The numerical solution $N$ satisfies the difference (20); replacing $u$ in (20) with $N$ in (3), one has

$$\frac{E_i^{m+1} + \varepsilon_i^{m+1} - E_i^m - \varepsilon_i^m}{\Delta t}$$
$$= -a\frac{E_{i+1}^m + \varepsilon_{i+1}^m - E_{i-1}^m - \varepsilon_{i-1}^m}{2\Delta x} \quad (21)$$

Similarly, replacing $u$ in (20) with $E$ gives

$$\frac{E_i^{m+1} - E_i^m}{\Delta t} = -a\frac{E_{i+1}^m - E_{i-1}^m}{2\Delta x}, \quad (22)$$

since $E$ clearly satisfies (20). Subtracting (22) from (21), one has

$$\frac{\varepsilon_i^{m+1} - \varepsilon_i^m}{\Delta t} = -a\frac{\varepsilon_{i+1}^m - \varepsilon_{i-1}^m}{2\Delta x}, \quad (23)$$

which proves that $\varepsilon$ also satisfies (20). The solution of (20) is stable if and only if the requirement in (7) is fulfilled.

Following the same argument as in the previous section, we are allowed to just handle one term of the Fourier series, (9), in order to probe the stability of (20) without loss in generality. Substituting (9) into (20) gives

$$\frac{e^{b(t+\Delta t)}e^{ik_n x} - e^{bt}e^{ik_n x}}{\Delta t}$$
$$= -a\frac{e^{bt}e^{ik_n(x+\Delta x)} - e^{bt}e^{ik_n(x-\Delta x)}}{2\Delta x} \quad (24)$$

Dividing (24) by $e^{bt}e^{ik_n x}$;

$$\frac{e^{b\Delta t} - 1}{\Delta t} = -a\frac{e^{ik_n\Delta x} - e^{-ik_n\Delta x}}{2\Delta x} \quad (25)$$

Combining (11) and (25);

$$\left|-\frac{1}{2}\frac{a\Delta t}{\Delta x}\left(e^{ik_n\Delta x} - e^{-ik_n\Delta x}\right) + 1\right| \leq 1 \quad (26)$$

Rewriting (26) in terms of $\xi_1$ and $\xi_2$, one has

$$\left|-\frac{1}{2}\frac{a\Delta t}{\Delta x}(\xi_1 - \xi_2) + 1\right| \leq 1 \quad (27)$$

The possible conditions are;

## A. Condition a

$$-\frac{1}{2}\frac{a\Delta t}{\Delta x}(\xi_1 - \xi_2) + 1 \leq 1 \quad (28)$$

Setting $\xi_1 = -1$, and $\xi_2 = 1$, (28) becomes

$$\frac{a\Delta t}{\Delta x} \leq 0 \quad (29)$$

## B. Condition b

$$\frac{1}{2}\frac{a\Delta t}{\Delta x}(\xi_1 - \xi_2) - 1 \leq 1 \quad (30)$$

Setting $\xi_1 = 1$, and $\xi_2 = -1$, (30) becomes

$$\frac{a\Delta t}{\Delta x} \leq 2 \quad (31)$$

Since both conditions, (29) and (31), must hold simultaneously, then not only the CFL criterion is not fulfilled but also (20) is found to be unconditionally unstable; the situation which justify the selection of Lax's difference scheme in the previous section as an option for the time derivative in (1). We will further examine whether such instability is solely due to the simple first-order forward difference scheme for the time derivative or the combination of schemes that represent both derivatives (i.e. time and spatial derivatives), in the following section.

## V. SIMPLE FIRST-ORDER DIFFERENCEs FOR THE TIME AND SPATIAL DERIVATIVES

In this section we consider two combinations involving simple *first-order forward* and *first-order rearward* difference schemes.

### A. Combination I

Replacing time and spatial derivatives with simple first-order finite differences to represent the derivatives in (1), we have, after some rearrangement

$$\frac{u_i^{m+1} - u_i^m}{\Delta t} = -a\frac{u_{i+1}^m - u_i^m}{\Delta x} \quad (32)$$

It can be shown that

$$\frac{\varepsilon_i^{m+1} - \varepsilon_i^m}{\Delta t} = -a\frac{\varepsilon_{i+1}^m - \varepsilon_i^m}{\Delta x}, \quad (33)$$

since (32) is satisfied by both $N$ and $E$. Substituting (9) into (33) and dividing the resulting equation by $e^{bt}e^{ik_n x}$ gives, after some rearrangement,

$$e^{b\Delta t} = -a\Delta t \frac{e^{ik_n \Delta x} - 1}{\Delta x} + 1 \quad (34)$$

The stability of the solution of (32) requires that

$$\left|\frac{-a\Delta t}{\Delta x}(\xi_1 - 1) + 1\right| \leq 1, \quad (35)$$

where $\xi_1 = e^{ik_n \Delta x}$.

#### 1) Condition a

$$\frac{-a\Delta t}{\Delta x}(\xi_1 - 1) \leq 0 \quad (36)$$

Setting $\xi_1 = -1$, (36) becomes

$$\frac{a\Delta t}{\Delta x} \leq 0 \quad (37)$$

#### 2) Condition b

$$\frac{a\Delta t}{\Delta x}(\xi_1 - 1) \leq 2 \quad (38)$$

Setting $\xi_1 = 1$ in (38), one has a trivial solution. Thus the solution of (32) is unconditionally unstable.

### B. Combination II

We now represent time and spatial derivatives in (1) with simple first-order *forward* and *rearward* finite-differences, respectively. The error $\varepsilon$ satisfies the difference equation;

$$\frac{u_i^{m+1} - u_i^m}{\Delta t} = -a\frac{u_i^m - u_{i-1}^m}{\Delta x} \quad (39)$$

It can be shown that

$$\frac{\varepsilon_i^{m+1} - \varepsilon_i^m}{\Delta t} = -a\frac{\varepsilon_i^m - \varepsilon_{i-1}^m}{\Delta x}, \quad (40)$$

and

$$\frac{e^{b(t+\Delta t)}e^{ik_n x} - e^{bt}e^{ik_n x}}{\Delta t} = -a\frac{e^{bt}e^{ik_n x} - e^{bt}e^{ik_n(x-\Delta x)}}{\Delta x} \quad (41)$$

Dividing the resulting difference equation for the round-off error $\varepsilon$ by $e^{bt}e^{ik_n x}$, the stability criterion is given by

$$e^{b\Delta t} = \left|-a\frac{\Delta t}{\Delta x}(1 - \xi_2) + 1\right| \leq 1, \quad (42)$$

where $\xi_2 = e^{-ik_n \Delta x}$.

#### 1) Condition a

$$a\frac{\Delta t}{\Delta x}(1-\xi_2) \geq 0 \quad (43)$$

Setting $\xi_2 = 1$ in (43) results in a trivial solution.

*2) Condition b*

$$a\frac{\Delta t}{\Delta x}(1-\xi_2) \leq 2 \quad (44)$$

Setting $\xi_2 = -1$ in (44);

$$a\frac{\Delta t}{\Delta x} \leq 1 \quad (45)$$

Thus the solution of (39) is conditionally stable, where the stability criterion is that of CFL.

## VI. CONCLUSION

We summarize all four cases of the finite-difference representations for the first-order hyperbolic wave equation in Table I.

TABLE I. DIFFERENCE REPRESENTATIONS OF THE WAVE EQUATION, BASED ON THE COMBINATIONS OF THE SCHEMES; LAX'S DIFFERENCE (LD), FORWARD DIFFERENCE (FD), REARWARD DIFFERENCE (RD), AND CENTRAL DIFFERENCE (CD)

| | Case | Time Derivative | Spatial Derivative | Remark |
|---|---|---|---|---|
| **Difference Scheme** | 1 | LD | CD | CFL criterion satisfied |
| | 2 | FD | CD | unconditionally unstable |
| | 3 | FD | FD | unconditionally unstable |
| | 4 | FD | RD | CFL criterion satisfied |

It has been shown in Case 1 and Case 2, that the first-order accurate LD is a suitable candidate for the FD replacement to represent the time derivative in the first-order wave equation while maintaining the CD for the spatial derivative, where the CFL is met. However, the FD seems to contribute to the stability of the corresponding difference equation solution if the spatial derivative is represented by the RD as shown in Case 4. This suggests the applicability of the FD in solving a hyperbolic partial differential equation.

ACKNOWLEDGMENT

The author is very grateful to the Department of Aeronautical Engineering, Faculty of Mechanical and Manufacturing Engineering, Universiti Tun Hussein Onn Malaysia where this work was initiated.